\documentclass[11pt]{article}
\usepackage[leqno]{amsmath}   
\usepackage{amsthm}  
\usepackage{amssymb}  
\usepackage{amscd}   

\newcommand{\proLog}{\operatorname{{\cal L}og}}

\newcommand{\Pic}{\operatorname{P}}

\newcommand{\ab}{\operatorname{ab}}

\newcommand{\Tr}{\operatorname{Tr}}

\newcommand{\Pol}{\operatorname{{\cal P}ol}}

\newcommand{\Gext}{\operatorname{{\cal P}ol}}

\newcommand{\Rh}{{\cal R}}
\newcommand{\Oh}{{\cal O}}

\newcommand{\I}{{\cal I}}
\newcommand{\1}{{\bf 1}}
\newcommand{\et}{{\operatorname{et}}}

\newcommand{\Lh}{{\cal L}}

\newcommand{\sbar}{\overline{s}}
\newcommand{\xbar}{\overline{x}}
\renewcommand{\tilde}{\widetilde}

\newcommand{\Norm}{\operatorname{N}}

\newcommand{\cl}{\operatorname{cl}}

\renewcommand{\Bbb}{\mathbb}  

\newcommand{\Q}{{\Bbb{Q}}}  
\newcommand{\Z}{{\Bbb{Z}}}  

\renewcommand{\epsilon}{\varepsilon}
\renewcommand{\rho}{\varrho}
\renewcommand{\bar}{\overline}

\newcommand{\isom}{\cong}       
\newcommand{\ohne}{\smallsetminus}
\newcommand{\Hom}{\operatorname{Hom}}
\newcommand{\Ext}{\operatorname{Ext}}
\newcommand{\prolim}{\varprojlim}

\newcommand{\bew}{\begin{proof}}
\newcommand{\bewende}{\end{proof}}
\newtheorem{lemma}{Lemma}[subsection]
\newtheorem{prop}[lemma]{Proposition}
\newtheorem{thm}[lemma]{Theorem}
\newtheorem{defn}[lemma]{Definition}
\newtheorem{Conditions}[lemma]{Conditions}

\newtheorem{cor}[lemma]{Corollary}

\newcommand{\verk}{{\scriptstyle\circ}}

\def\ab{\operatorname{{ab}}}
\title{\vskip -25mm A note on polylogarithms on curves and abelian schemes}   
\author{Guido Kings}
\date{}
\begin{document}

\maketitle
\setcounter{section}{0}

\section*{Introduction}

Cohomology classes defined by polylogarithms have been one of the main tools to study special
values of $L$-functions. Most notably, they play a decisive role in the study of the
Tamagawa number conjecture for abelian number fields (\cite{Be2}, \cite{Del1}, \cite{HW} \cite{Hu-Ki}), 
CM elliptic curves (\cite{Deninger}, \cite{Ki2}) and modular forms
(\cite{Be1}, \cite{Kato}).

Polylogarithms have been defined for relative curves by Beilinson and 
Levin (unpublished) and for abelian schemes by Wildeshaus \cite{Wi} in the context of
mixed Shimura varieties. In general, the nature of these extension classes is not well understood.
 
The aim of this note is to show that there is a close connection
between the polylogarithm extension on curves and on abelian schemes. 
It turns out that the polylog on an abelian scheme is roughly the
push-forward of the polylog on a sub-curve. If we apply this to
the embedding of a curve into its Jacobian, we can give a more
precise statement: the polylog on the Jacobian is the cup product 
of the polylog on the curve with the fundamental class of the curve
(see theorem \ref{cupprodthm}). 
With this result it is possible to understand the nature of 
the polylog extension on abelian schemes in a better way. 

The polylog extension on 
curves has the advantage of being a one extension of lisse 
sheaves. Thus, itself can be represented by a lisse sheaf. The 
polylog extension on the abelian scheme on the contrary is a
$2d-1$ extension, where $d$ is the relative dimension
of the abelian scheme.
%


The contents of this note is as follows: To simplify the
exposition we only treat the \'etale realization. First we define the polylog extension on
curves and abelian schemes in a unified way for integral coefficients. To our knowledge this
and the construction on curves is not published but goes back to an earlier version of \cite{Be-Le}. 
The case of abelian schemes is treated in \cite{Wi} (for $\Q_l$-sheaves),
which we mildly generalize to $\Z/l^r\Z$- and $\Z_l$-sheaves. All the main ideas are
of course already in \cite{Be-Le}.

The second part gives three important properties of the polylog extension,
namely compatibility with base change, norm compatibility and the splitting principle.

In the last part we show that the push-forward of the polylog on a sub-curve of an abelian scheme
gives the polylogarithmic extension on the abelian scheme and prove our main theorema about the
polylog on the Jacobian.

%

\section{Definition of the polylogarithm extension}

The first part of this  paper recalls the definition of the
polylogarithmic extension for curves and abelian schemes. 

The case of elliptic curves was treated by Beilinson and Levin \cite{Be-Le}
in analogy with the cyclotomic case considered by Beilinson and Deligne.
An earlier version of \cite{Be-Le} contained also the case of general curves.
Polylogarithmic extensions
for abelian schemes and more generally certain semi-abelian schemes were
first considered by Wildes\-haus in \cite{Wi} in the context of mixed 
Shimura varieties.

\subsection{The logarithm sheaf}\label{seclogarithm}
In this section we recall the definition of  the logarithm
sheaf for curves and abelian schemes. 

Let $S$ be a connected scheme, and $l$ be a prime number invertible on $S$.
We fix a base ring $\Lambda$ which is either $\Z/l^r\Z$ or $\Z_l$. 
The cohomology in this paper is always continuous cohomology in the
sense of Jannsen \cite{Ja}.

\begin{defn}\label{curvedef} A {\em curve } is  a  smooth 
proper morphism $\pi :C\to S$ together with a section $e:S\to C$,
such that the geometric fibers $C_{\sbar}$ of 
$\pi$ are  connected curves of genus $\ge 1$. 
\end{defn}
In addition to the curves we consider
also abelian schemes $\pi:A\to S$ with unit section $e:S\to A$.
 For brevity we use the following
notation: $\pi:X\to S$  will denote either a curve in the sense of 
\ref{curvedef} or an abelian scheme over $S$. The  section will be denoted 
by $e:S\to X$. The relative dimension of $X/S$ is $d$. 

Let us describe a $\Lambda$-version of the theory in \cite{Wi} I, chapter 3.
In the case of an elliptic curve this coincides with \cite{Be-Le}.
Let $\bar{s}$ be a geometric point of $S$ 
and denote by $\bar{x}:=e(\bar{s})$ the corresponding geometric
point of $X$. Denote the fiber over $\bar{s}$ by $X_{\bar{s}}$
and consider the split exact sequence of fundamental groups
\begin{equation}\label{relfundgroup}
1\to \pi_1'(X_{\bar{s}},\bar{x})\to
\pi_1'(X,\bar{x})\xrightarrow{\pi_*}\pi_1(S,\bar{s})\to 1
\end{equation}
(cf. \cite{SGA1} XIII 4.3),
where $ \pi_1'(X_{\bar{s}},\bar{x})$ is the largest pro-$l$-quotient 
of $ \pi_1(X_{\bar{s}},\bar{x})$ and if $\ker(\pi_*)/N$
denotes the largest pro-$l$-quotient of $\ker(\pi_*)$,
then $\pi_1'(X,\bar{x}):=\pi_1(X,\bar{x})/N$.
The splitting is given by $e_*$. Now 
$  \pi_1'(X_{\bar{s}},\bar{x})$ is a pro-finite group and 
we fix a fundamental system of open neighborhoods $\Gamma_{j}$ of the identity with $j\in J$,
 such that
\[
 \pi_1'(X_{\bar{s}},\bar{x})=\prolim_j \pi_1'(X_{\bar{s}},\bar{x})/\Gamma_j.
\]
Define 
\[
H_j:=\pi_1'(X_{\bar{s}},\bar{x})/\Gamma_j
\]
and let us agree that in the case where $X$ is an abelian scheme we choose
the $\Gamma_j$ in such a way that $H_j=\ker[l^j]$ is the kernel of the
$[l^j]$-multiplication.
Let us also fix a projective system $X_j$ of \'etale $H_j$-torsors
(i.e. Galois coverings of $X$ with group $H_j$)
such that
\[
\begin{CD}
H_j@>h_j>>X_j\\
@VVV@VVp_jV\\
S@>e>>X,
\end{CD}
\]
is Cartesian. 
In the case of an abelian scheme we take  $X_j=A$ and $p_j=[l^j]$.
For $j'\to j$ we have the trace map
\[
{p_{j'}}_*\Lambda\to {p_j}_*\Lambda
\]
and we define:
\begin{defn}
The {\em logarithm  sheaf} is the \'etale sheaf
\[
\proLog_{X,\Lambda}:=\prolim_j({p_j}_*\Lambda)
\]
where the transition maps are the 
above trace maps.
\end{defn}
The stalk of $\proLog_{X,\Lambda}$ at $\bar{x}$ is just the 
Iwasawa algebra of the profinite group  $\pi_1'(X_{\bar{s}},\bar{x})$
$$
	\proLog_{X,\Lambda,\bar{x}}=\Lambda[[ \pi_1'(X_{\bar{s}},\bar{x})]]
$$ 
with the canonical action of the semi-direct product $\pi_1'(X,\bar{x})$ given by
multiplication on $\pi_1'(X_{\bar{s}},\bar{x})$ and by conjugation on the quotient
$\pi_1(S,\bar{s})$.

Let 
\[
\Rh_{X,\Lambda}:=e^*\proLog_{X,\Lambda}
\]
be the pull-back of $\proLog_{X,\Lambda}$ along the unit section $e$. 

Note that that there is a canonical map $\1:\Lambda\to \Rh_{X,\Lambda}$
and that $\Rh_{X,\Lambda}$ has a ring structure given by group multiplication 
as usual. We denote by 
\[
\I_{X,\Lambda}:=\ker(\Rh_{X,\Lambda}\to \Lambda)
\]
the augmentation ideal.
The logarithm sheaf has a canonical action of $\pi^*\Rh_{X,\Lambda}$
(again induced by group multiplication)
\[
\pi^*\Rh_{X,\Lambda}\otimes_{\Lambda}\proLog_{X,\Lambda}\to \proLog_{X,\Lambda},
\]
which defines on $\proLog_{X,\Lambda}$ the structure of an $\pi^*\Rh_{X,\Lambda}$-torsor. 

It is very useful to consider the abelianized version of the logarithm sheaf.
\begin{defn}\label{ablogdefn}
	Let $\pi_1'(X_{\bar{s}},\bar{x})^{\ab}$ be the maximal abelian quotient
	of $\pi_1'(X_{\bar{s}},\bar{x})$ and define the \emph{abelian logarithm sheaf}
	to be the lisse sheaf defined by the $\pi_1'(X,\bar{x})$-representation
	$$
		\proLog_{X,\Lambda,\bar{x}}^{\ab}:= \Lambda[[ \pi_1'(X_{\bar{s}},\bar{x})^{\ab}]].
	$$
\end{defn}

Note that in the our case $\pi_1'(X_{\bar{s}},\bar{x})^{\ab}\isom \Z_l^r$ is a free modul over $\Z_l$ of
rank twice the genus of the curve or twice the dimension of the ablian scheme.
In particular, 
$$
	\Lambda[[ \pi_1'(X_{\bar{s}},\bar{x})^{\ab}]]\isom \Lambda[[x_1,\ldots,x_r]]
$$
 is isomorphic to a power series ring in $r$ variables.

The whole theory of the polylogarithm sheaves relies on the fact that 
the higher direct images of $\proLog_{X,\Lambda}$ can be computed for curves and 
abelian schemes.
\begin{thm}[\cite{Be-Le},\cite{Be-Le2},\cite{Wi}]\label{highergenim}
Let $d$ be the relative dimension of $X=C,A$ over $S$.
Then the \'etale sheaf 
\[
R^i\pi_*\proLog_{X,\Lambda}
\] 
is zero for $i\neq 2d$ and 
\[
R^{2d}\pi_*\proLog_{X,\Lambda}\isom \Lambda(-d).
\] 
\end{thm}
\bew We use the spectral sequence \cite{Ja} 3.10. We have to compute the transition maps in the inverse system 
$(R^i\pi_*{p_j}_*\Lambda)_j$, which by Poincar\'e duality can by written as an inductive system
\[
R^{2d-i}\pi_!{p_j}_*\Lambda(d).
\]
Here the transition maps are induced by the pull-back maps $j'\to j$
\[
{p_j}_*\Lambda(d)\to {p_{j'}}_*\Lambda(d).
\]
Let us first show that the inductive system 
$(R^{2d-i}\pi_!{p_j}_*\Lambda(d))_j$  is zero
for $i< 2d$. By base change we may assume that $S$ is the spectrum of an
algebraically closed  
field, so that we have to compute the transition maps for
\[
H^{2d-i}_{\et}(X_j, \Lambda(d)).
\]
We have 
\[
H^{1}_{\et}(X_j, \Lambda(d))=\Hom(\pi_1(X_j),\Lambda(d))=\Hom(\Gamma_j,\Lambda(d)).
\]
For every homomorphism in $\Hom(\Gamma_j,\Lambda(d))$ there is an $j'$, such
that its restriction to $\Gamma_{j'}$ is trivial. Thus the maps 
\[
 H^{1}_{\et}(X_j, \Lambda(d))\to H^{1}_{\et}(X_{j'}, \Lambda(d))
\]
are zero.
For an abelian scheme we have
\[
H^{2d-i}_{\et}(X_j, \Lambda(d))\isom \Lambda^{2d-i}H^{1}_{\et}(X_j, \Lambda(d))
\]
hence we get the same result for 
\[
H^{2d-i}_{\et}(X_j, \Lambda(d))=0
\] 
for $i<2d$. Now 
$ H^{0}_{\et}(X_j, \Lambda(d))=\Lambda(d)$ is constant,  so that 
the natural map $R^{2d}\pi_!{p_j}_*\Lambda(d)\to \Lambda$ is an isomorphism.
Thus the result is proven  for abelian schemes. 
If $X$ is a curve it remains to consider
\[
H^{2}_{\et}(X_j, \Lambda(1))
\]
We have $H^{2}_{\et}(X_j, \Lambda(1))=\Lambda$ and the transition maps
are given by multiplication with the degree of $X_{j'}\to X_j$. Thus for
$j'$ large enough the transition maps are zero.
\bewende

\subsection{The polylogarithm extension}

Let $X$ be a curve or an abelian scheme and $U:=X\ohne e(S)$.
Denote by $\pi_U:U\to S$ the restriction of
$\pi $ to $U$.
We let 
\[
j:U\to X
\]
be the open immersion of $U$ into $X$. The restriction 
of $\proLog_{X,\Lambda}$ to 
$U$ is denoted  by $\proLog_{U,\Lambda}$. 

\begin{prop} Let $d$ be the relative dimension of $X/S$ and recall that
$\I_{X,\Lambda}$ is the kernel of the augmentation map of $\Rh_{X,\Lambda}$.
Then the \'etale sheaf
\[
R^i\pi_{U*}\proLog_{U,\Lambda}(d)
\] 
is zero for $i\neq 2d-1$ and for $i=2d-1$ there is 
an isomorphism of sheaves
\[
R^{2d-1}\pi_{U*}\proLog_{U,\Lambda}(d)\isom \I_{X,\Lambda}.
\]
\end{prop}
\bew
This follows immediately from the localization sequence 
\[
R^{i-1}\pi_{U*}\proLog_{U,\Lambda}(d)\to R^{i}e^!\proLog_{X,\Lambda}(d)\to 
R^{i}\pi_{*}\proLog_{X,\Lambda}(d)
\]
and the purity isomorphism $e^!\proLog_{X,\Lambda}(d)\isom e^*\proLog_{X,\Lambda}[-2d]
=\Rh_{X,\Lambda}[-2d]$.
By theorem \ref{highergenim} the sheaf $ R^{i}\pi_{*}\proLog_{X,\Lambda}(d)$ is zero
 for $i\neq 2d$, so that we get an exact sequence
\[
 R^{2d-1}\pi_{*}\proLog_{X,\Lambda}(d)\to  R^{2d-1}\pi_{U*}\proLog_{U,\Lambda}(d)\to\Rh_{X,\Lambda}\to 
R^{2d}\pi_{*}\proLog_{X,\Lambda}(d)\to 0.
\]
The identification $ R^{2d}\pi_{*}\proLog_{X,\Lambda}(d)\isom \Lambda$
gives an identification of the last map with the augmentation map
$ \Rh_{X,\Lambda}\to \Lambda$ and $R^{2d-1}\pi_{U*}\proLog_{U,\Lambda}(d)$
becomes isomorphic to the augmentation ideal.
\bewende
Consider the extension 
\[
\Ext^{i}_{U}(\pi_U^*\I_{X,\Lambda},\proLog_{U,\Lambda}(d)).
\]
\begin{cor}\label{extiso}
There is an isomorphism
\[
\Ext^{2d-1}_{U}(\pi_U^*\I_{X,\Lambda},\proLog_{U,\Lambda}(d))\isom \Hom_S(\I_{X,\Lambda},\I_{X,\Lambda}) 
\]
given by the edge morphism in the Leray spectral sequence for 
$R\pi_{U*}$. 
\end{cor}
\bew 
In  the Leray spectral sequence for $ R\pi_{U*}$ all higher
direct images except $R^{2d-1}\pi_{U*}\Lh_j|_U(d)$ are 
zero. Thus all higher $\Ext$ groups vanish  except 
\[
\Hom_S(\I_{X,\Lambda},R^{2d-1}\pi_{U*}\proLog_{U,\Lambda}(d)).
\]
The above theorem identifies $R^{2d-1}\pi_{U*}\proLog_{U,\Lambda}(d)$
with $\I_{X,\Lambda}$ and the result follows.
\bewende
\begin{defn}\label{poldefn}
The {\em large polylogarithm  extension } on $X$ is the  extension
class $\Pol_{X,\Lambda}$ in 
\[
\Ext^{2d-1}_{U}(\pi_U^*\I_{X,\Lambda}, \proLog_{U,\Lambda}(d))\isom \Hom_S(\I_{X,\Lambda},\I_{X,\Lambda})
\]
corresponding to the identity in $ \Hom_S(\I_{X,\Lambda},\I_{X,\Lambda})$.
\end{defn}
\begin{rem} a) In the case of an elliptic curve the polylogarithm
is the one considered by Beilinson and Levin \cite{Be-Le}.
In the case of an abelian scheme our definition gives a 
$\Z/l^r\Z$-version of the construction in Wildeshaus \cite{Wi}.
\\
b) This polylog extension should be more precisely called the 
\'etale realization of the polylog. 
The Hodge-realization can be defined in a similar way.\\
c) In the case of an abelian
scheme, this class is in the
image of the regulator coming from K-theory
(see \cite{Ki}). 
\end{rem}

It is useful to make also the following definition:

\begin{defn}\label{abpoldefn}
The {\em abelian polylogarithm  extension } on $X$ is the  extension
class $\Pol_{X,\Lambda}^{\ab}$ in 
\[
\Ext^{2d-1}_{U}(\pi_U^*\I_{X,\Lambda}, \proLog_{U,\Lambda}^{\ab}(d)),
\]
which is the image of $\Pol_{X,\Lambda}$ under the canonical map
$$
	\Ext^{2d-1}_{U}(\pi_U^*\I_{X,\Lambda}, \proLog_{U,\Lambda}(d))\to \Ext^{2d-1}_{U}(\pi_U^*\I_{X,\Lambda}, \proLog_{U,\Lambda}^{\ab}(d)).
$$
\end{defn}

\section{Properties of the polylog}

We consider compatibilities of the polylog, namely its behavior
under base change and finite \'etale morphisms. Finally, we explain
the splitting principle in our situation.

\subsection{Compatibility with base change}
Let $S'\to S$ be a scheme over $S$ and 
$X'$ be the fiber product $X\times_SS' $, where $X$ is $C$ or
$A$ as usual.
We get a Cartesian diagram
\[
\begin{CD}X'@>f>>X\\
@VV\pi'V @VV\pi V \\
S'@>g>>S.
\end{CD}
\]
and let $U':=X'\ohne e'(S')$ and $U:=X\ohne e(S)$.
Here $e$ and $e'$ are the unit sections of $X$ and $X'$ respectively.
The obvious map of fundamental groups
$\pi_1'(X_{\sbar},\xbar)\to \pi'_1(X'_{\sbar'},\xbar')$
induces
\[
\proLog_{X'}\to f^*\proLog_X
\]
which is  called the canonical map.
\begin{lemma}\label{BClog}  The canonical map
\[
\proLog_{X'}\to f^*\proLog_X
\]
is an isomorphism. 
\end{lemma}
\bew 
It suffices to remark that pull-back by $f$ induces an
isomorphism of the relative  fundamental groups 
\[
\pi_1'(X_{\sbar},\xbar)\isom \pi'_1(X'_{\sbar'},\xbar')
\]
from equation (\ref{relfundgroup}) in section \ref{seclogarithm}. Here we
use of course the base points
$f(\xbar')=\xbar$ and $f(\sbar')=\sbar$.
\bewende
In particular the canonical map induces an isomorphism
\[
{\pi'}^*\I_{X'}\isom f^*\pi^*\I_X
\]
so that we have
\[
\Ext^{2d-1}_U(f^*\pi_U^*\I_X,f^*\proLog_U(d))\isom
\Ext^{2d-1}_{U'}({\pi'_{U'}}^*\I_{X'},\proLog_{U'}(d)).
\]
\begin{cor} Via this identification of $\Ext$-groups
\[
f^*\Gext_X=\Gext_{X'}.
\]
\end{cor}
\bew We have a commutative diagram
\[
\begin{CD}
\Ext^{2d-1}_U(\pi_U^*\I_X,\proLog_U(d))@>f^*>>
\Ext^{2d-1}_{U'}({\pi'_{U'}}^*\I_{X'},\proLog_{U'}(d))\\
@VV\isom V@VV\isom V\\
\Hom_S(\I_X, \I_X)@>{g}^*>>\Hom_{S'}(\I_{X'},\I_{X'})
\end{CD}
\]
and the identity is mapped to the identity under 
${g}^*$. 
\bewende

\subsection{Norm compatibility for finite \'etale morphisms}
Let $f:X'\to X$ be a finite \'etale pointed morphism  over $S$, 
i.e. $f\verk e'=e$, and denote by $\pi'$ and $\pi$ the
structure maps of $X'$ and $X$. Let $U$ and $U'$ be as
above. Define $Z'$ and $\widetilde{U}$ by the Cartesian diagram:
\[
\begin{CD}
Z'@>\epsilon>>X'@<\tilde{j}<<\widetilde{U}\\
@VVfV@VVfV@VVfV\\
S@>e>>X@<<<U.
\end{CD}
\]
Observe that $\widetilde{U}\subset U'$ and denote by  $U'\xrightarrow{j'}X'$
the open immersion.
Restriction to $\widetilde{U}$ gives a map
\begin{equation}\label{eq1}
\Ext^{2d-1}_{U'}({\pi'_{U'}}^*\I_{X'},\proLog_{U'}(d))\to 
\Ext^{2d-1}_{\widetilde{U}}
({\pi'_{\widetilde{U}}}^*\I_{X'},\proLog_{\widetilde{U}}(d)).
\end{equation}
On the other hand we have an adjunction map (${\pi'_{\widetilde{U}}}^*=
f^*\verk {\pi_{{U}}}^*$)
\begin{equation}\label{eq2}
\Ext^{2d-1}_{\widetilde{U}}(f^*{\pi_{{U}}}^*\I_{X'},\proLog_{\widetilde{U}}(d))
\to 
\Ext^{2d-1}_{U}({\pi_U}^*\I_{X'},f_*\proLog_{\widetilde{U}}(d)).
\end{equation}
If we apply $f_*$ to the canonical map $\proLog_{X'}\to f^*\proLog_X$
we get a map
\[
f_*\proLog_{\widetilde{U}}\to f_*f^*\proLog_U\xrightarrow{\Tr}\proLog_U
\]
and hence a map 
\begin{equation}\label{eq3}
\Ext^{2d-1}_{U}({\pi_U}^*\I_{X'},f_*\proLog_{\widetilde{U}}(d))\to
\Ext^{2d-1}_U(\pi_U^*\I_{X'},\proLog_U(d)).
\end{equation}
Denote by 
$\Norm_f$ the resulting composition of (\ref{eq1}), (\ref{eq2}) and
(\ref{eq3}) 
\[
\Norm_f:\Ext^{2d-1}_{U'}({\pi'_{U'}}^*\I_{X'},\proLog_{U'}(d))\to 
\Ext^{2d-1}_U(\pi_U^*\I_{X'},\proLog_U(d)).
\]
Consider also the map $\I_{X'}\to \I_X$
induced by the canonical map $\proLog_{X'}\to f^*\proLog_X$ and let
\[
\alpha^*:\Ext^{2d-1}_U(\pi_U^*\I_{X},\proLog_U(d))\to
\Ext^{2d-1}_U(\pi_U^*\I_{X'},\proLog_U(d))
\]
be the associated pull-back.
\begin{prop}{\em (Norm compatibility)} With the maps $\Norm_f$ and $\alpha^*$ defined 
above
\[
\Norm_f(\Gext_{X'})=\alpha^*(\Gext_X).
\]
\end{prop}
\bew First of all we remark that
\[
R^i\pi_{\widetilde{U},*}\proLog_{\widetilde{U}}
\]
is zero for $i\neq 2d-1$ and isomorphic
to 
\[
\ker(f_*\epsilon^*\proLog_{X'}\to \Lambda)
\]
 for $\epsilon:Z'\hookrightarrow X'$ and $i=2d-1$.
This follows as in corollary \ref{extiso} from the localization sequence
\[
\epsilon_*\epsilon^!\proLog_{X'}(d)\to \proLog_{X'}(d)\to R\tilde{j}_*\proLog_X(d)
\]
and the purity isomorphism 
$\epsilon^!\proLog_{X'}(d)\isom \epsilon^*\proLog_{X'}[-2d]$.
Thus the map $\Norm_f$ on $\Ext$-groups can be identified with the following
composition of 
maps on $\Hom$-groups:
\[
\Hom_S(\I_{X'},\I_{X'})\to 
\Hom_S(\I_{X'},\ker(f_*\epsilon^*\proLog_{X'}\to \Lambda))
\to \Hom_S(\I_{X'},\I_{X}),
\]
where the first map is induced by the natural inclusion of 
$\I_{X'}$ into the kernel of $f_*\epsilon^*\proLog_{X'}\to \Lambda$ and
the second is induced by the trace map.
We have to show that the identity in $\Hom_S(\I_{X'},\I_{X'})$
maps to the canonical map in $\Hom_S(\I_{X'},\I_{X})$ or that
\[
{e'}^*\proLog_{X'}\to f_*\epsilon^*\proLog_{X'}\to e^*\proLog_X
\]
is the canonical map. As the first map is $f_*$ applied to
$e'_*{e'}^*\proLog_{X'}\to \epsilon_*\epsilon^*\proLog_{X'}$ this is clear. 
\bewende
\subsection{The splitting principle}

Let $\phi:X\to X'$ be a finite Galois covering over $S$ (we assume that $\phi$ is a pointed map)
with covering group $G$, 
where $X,X'$ are either curves or abelian schemes.
Then we have a map
$$
	\phi_*:\proLog_{X,\Lambda}\to \phi^*\proLog_{X',\Lambda}
$$
induced by the canonical map $\phi_*:\pi_1'(X_{\bar{s}},\bar{x})\to \pi_1'(X'_{\bar{s}},\bar{x})$
of relative fundamental groups.

\begin{lemma}[Splitting principle]
	Assume in the above situation that the order of $G$ is prime to $l$, then
	$$
		\phi_*:\proLog_{X,\Lambda}\to \phi^*\proLog_{X',\Lambda}
	$$
	is an isomorphism. In particular, for any $x\in X$, which is in the $G$-orbit of
	$e$, one has a canonical isomorphism
	$$
		\Rh_{X,\Lambda}\isom x^*\proLog_{X,\Lambda}.
	$$
\end{lemma}

\bew
Obvious, as in this case $\phi_*:\pi_1'(X_{\bar{s}},\bar{x})\to \pi_1'(X'_{\bar{s}},\bar{x})$
is an isomorphism.
\bewende
This lemma is very useful in the case of the abelian polylog. One gets
$$
	x^*\Pol_{X,\Lambda}^{\ab}\in \Ext^{2d-1}_{U}(\pi_U^*\I_{X,\Lambda}, \Rh_{X,\Lambda}^{\ab}(d))
$$
and $\Rh_{X,\Lambda}^{\ab}$ is isomorphic to a power series ring over $\Lambda$.

\section{The connection between the abelian polylogarithm on curves and on
abelian schemes}
In this section we show that the polylogarithm on an abelian
scheme is induced from the abelian polylogarithm on a sub-curve.

\subsection{The polylog on an abelian scheme as push-forward}
We consider curves with the following:
\begin{Conditions}\label{cond}
\begin{itemize}
\item[a)] $i:C\hookrightarrow A$ is a closed embedding and $C$
is is  a curve as in definition \ref{curvedef}. 
\item[b)] The map $\rho:R^1\pi_{C*}\Lambda\to R^1\pi_{A*}\Lambda$ 
induced by $i$ is surjective
\item[c)] The section $e$ of $C$ is mapped under $i$ 
to the unit section $e$ of $A$
\end{itemize}
\end{Conditions}
The structure maps of $C$ and $A$ will be denoted by
$\pi_C$ and $\pi_A$.

\begin{rem}
The typical case to consider here is  $A=\Pic^0(C/S)$ and
$i:C\to A$ the embedding $t\mapsto \Oh((e)-(t))$.
\end{rem}

\begin{lemma} Let $i:C\hookrightarrow A$ satisfy the conditions \ref{cond},
then the canonical map 
\[
\rho:\proLog_C\to i^*\proLog_A
\]
is surjective (and factors through $\proLog_C^{\ab}$). In particular, the canonical maps $\rho:\Rh_C\to \Rh_A$ 
and $\rho:\I_C\to \I_A$ are surjective.
\end{lemma}
\bew The condition \ref{cond} b) implies that the pull-back of
the coverings $A_j\to A$ to $C$ are  quotients of the maximal
pro-$l$-covering of $C$. Applying this to the definition of 
$\proLog_C$ and $\proLog_A$ gives the desired result.
\bewende

Consider the adjunction map
\[
i_*i^!\proLog_A\to\proLog_A
\] 
and observe that by purity 
$i^!=i^*[-2d+2](-d+1)$ as $C$ is smooth, so that the
canonical map $\rho$ gives
\[
\lambda:i_*\proLog_C\to\proLog_A[2d-2](d-1).
\]
This gives us two diagrams
\[\begin{CD}
\Ext^1_{C\ohne e(S)}(\pi_C^*\I_C,\proLog_C(1))@>{\alpha}>>
\Ext^{2d-1}_{A\ohne e(S)}(\pi_A^*\I_C,\proLog_A(d))\\
@VV\isom V@VV\isom V\\
\Hom_S(\I_C,\I_C)@>\rho_*>>
\Hom_S(\I_C,\I_A),
\end{CD}
\]
where $\alpha$ is given by the adjunction map composed with $\lambda$
and 
\[\begin{CD}
\Ext^{2d-1}_{A\ohne e(S)}(\pi_A^*\I_A,\proLog_A(d))@>{\beta}>>
\Ext^{2d-1}_{A\ohne e(S)}(\pi_A^*\I_C,\proLog_A(d))\\
@VV\isom V@VV\isom V\\
\Hom_S(\I_A,\I_A)@>\rho^* >>
\Hom_S(\I_C,\I_A),
\end{CD}
\]
where $\beta$ is the pull-back by $\rho$. Note that the isomorphisms
for the right vertical arrows follow from a trivial extension of 
corollary \ref{extiso}.
\begin{prop}\label{polcpola} With the above notation, 
\[
\alpha(\Pol_C)=\beta(\Pol_A).
\]
\end{prop}
\bew
Results immediately from the above diagrams.
\bewende

\subsection{ The polylog on the Jacobian  as cup-product}
The push forward of the polylog on a curve to its Jacobian
can also be written as a cup-product. This sheds some light on the nature
of the extensions involved. Note that the polylog on curves
is a one extension of lisse sheaves, hence itself represented
by a lisse sheaf. 
It turns out, that the polylog on the Jacobian is the 
Yoneda product of this extension with the fundamental class of the curve.

Let $J=\Pic^0(C/S)$ be the Jacobian of $C$ and $i:C\to J$ be an embedding
satisfying the conditions in \ref{cond}. Let $g$ be the genus of $C$.

Consider in 
\[
\Ext^{2g-2}_{J\ohne e(S)}(i_*i^*\Lambda, \Lambda(g-1))\isom 
\Hom_{C\ohne e(S)}(\Lambda, \Lambda)
\]
the class $cl(C)$ corresponding to the identity.
This is the fundamental class of $C\ohne e(S)$ in $J\ohne e(S)$. If we tensor
this class with $\proLog_J(1)$ we get
\[
cl(C)\otimes\proLog_J(1)\in \Ext^{2g-2}_{J\ohne e(S)}(i_*i^*\proLog_J(1), 
\proLog_J(g)).
\]
On the other hand we have 
\[
i_*\Pol_C\in \Ext^1_{J\ohne e(S)}(i_*\pi_C^*\I_C, i_*\proLog_C(1))
\]
and with the canonical map $\rho:i_*\proLog_C(1)\to i_*i^*\proLog_J(1)$
we get 
\[
\rho_*i_*\Pol_C\in \Ext^{1}_{J\ohne e(S)}(i_*\pi_C^*\I_C,i_*i^*\proLog_J(1)).
\]
Denote by $\gamma(\Pol_C)$ the image of this class under the pull-back
by the adjunction
map $\pi_J^*\I_C\to i_*\pi_C^*\I_C$ in
\[
\Ext^{1}_{J\ohne e(S)}(\pi_J^*\I_C,i_*i^*\proLog_J(1)).
\]
The Yoneda product gives a class 
\[
\gamma(\Pol_C)\cup (cl(C)\otimes\proLog_J(1))\in
\Ext^{2g-1}_{J\ohne e(S)}(\pi_J^*\I_C, 
\proLog_J(g)).
\]
The canonical 
map $\rho:\I_C\to \I_J$ induces a map $\pi_J^*\I_C\to \pi_J^*\I_J$,
which gives
\[
\Ext^{2g-1}_{J\ohne e(S)}(\pi_J^*\I_J, 
\proLog_J(g))\xrightarrow{\beta}\Ext^{2g-1}_{J\ohne e(S)}(\pi_J^*\I_C, 
\proLog_J(g)).
\]

\begin{thm}\label{cupprodthm} In $\Ext^{2g-1}_{J\ohne e(S)}(\pi_J^*\I_C, \proLog_J(g))$
holds the equality
\[
\beta(\Pol_J)=\gamma(\Pol_C)\cup (cl(C)\otimes\proLog_J(1)).
\]
\end{thm}
\bew 
Let us write the extension classes as morphisms in the derived category.
We have 
\begin{align*}
\gamma(\Pol_C):&\pi_J^*\I_C\longrightarrow i_*i^*\proLog_J(1)[1]\\
\cl(C)\otimes\proLog_J(1)[1]:&i_*i^*\proLog_J(1)[1]\longrightarrow
\proLog_J(g)[2g-1]
\end{align*}
so that
\[
\gamma(\Pol_C)\cup (cl(C)\otimes\proLog_J(1)[1]):i_*\pi_C^*\I_C\longrightarrow
\proLog_J(g)[2g-1]
\]
is the composition. 
This morphism is the adjoint of 
\[
\pi_C^*\I_C\xrightarrow{\Pol_C} \proLog_C(1)[1]\to i^*\proLog_J(1)[1]
\]
because the adjoint of $cl(C)\otimes\proLog_J(1)[1]$ is the identity.
We have a commutative diagram
\[
\begin{CD}
\pi_J^*\I_C@>>> i_*i^*\proLog_J(1)[1]\\
@VVV@VVV\\
i_*\pi_C^*\I_C@>>> \proLog_J(g)[2g-1]
\end{CD}
\]
where the left and right vertical arrows are adjunction maps and
the horizontal arrows are adjoint to $\Pol_C\verk\rho$. 
The composition of these
maps is the morphism in the proposition and 
\[
\pi_J^* \I_C\to i_*i^*\proLog_J(1)[1]\to \proLog_J(2g-1)[2g-1]
\]
is just $\alpha(\Pol_C)$ with the notation of proposition \ref{polcpola}.
Our result follows from this proposition.
\bewende

\noindent Guido Kings\\
NWF-I Mathematik \\
Universit\"at Regensburg\\
93040 Regensburg\\
guido.kings@mathematik.uni-regensburg.de

\end{document}